\documentclass[journal]{IEEEtran}

\usepackage[cmex10]{amsmath}
\usepackage{amsfonts, amssymb, graphicx, cite}

\newtheorem{theorem}{Theorem}
\newtheorem{remark}{Remark}
\DeclareMathOperator{\sat}{sat}

\begin{document}

\title{Voltage/Pitch Control for Maximization and\\ Regulation of Active/Reactive Powers in\\ Wind Turbines with Uncertainties}

\author{Yi Guo, S. Hossein Hosseini, John N. Jiang, Choon Yik Tang, and Rama G. Ramakumar
\thanks{Yi Guo, S. Hossein Hosseini, John N. Jiang, and Choon Yik Tang are with the School of Electrical and Computer Engineering, University of Oklahoma, Norman, OK 73019 USA (e-mail: \{yi.guo,s.h.hosseini,jnjiang,cytang\}@ou.edu).}
\thanks{Rama G. Ramakumar is with the School of Electrical and Computer Engineering, Oklahoma State University, Stillwater, OK 74078 USA (e-mail: ramakum@okstate.edu).}
\thanks{This work was supported by the National Science Foundation under grants ECCS-0926038 and ECCS-0955265.}
\thanks{This paper is a preprint of a paper submitted to IET Renewable Power Generation and is subject to Institution of Engineering and Technology Copyright. If accepted, the copy of record will be available at IET Digital Library.}}

\maketitle

\begin{abstract}
This paper addresses the problem of controlling a variable-speed wind turbine with a Doubly Fed Induction Generator (DFIG), modeled as an electromechanically-coupled nonlinear system with rotor voltages and blade pitch angle as its inputs, active and reactive powers as its outputs, and most of the aerodynamic and mechanical parameters as its uncertainties. Using a blend of linear and nonlinear control strategies (including feedback linearization, pole placement, uncertainty estimation, and gradient-based potential function minimization) as well as time-scale separation in the dynamics, we develop a controller that is capable of maximizing the active power in the {\em Maximum Power Tracking} (MPT) mode, regulating the active power in the {\em Power Regulation} (PR) mode, seamlessly switching between the two modes, and simultaneously adjusting the reactive power to achieve a desired power factor. The controller consists of four cascaded components, uses realistic feedback signals, and operates without knowledge of the $C_p$-surface, air density, friction coefficient, and wind speed. Finally, we show the effectiveness of the controller via simulation with a realistic wind profile.
\end{abstract}


\section{Introduction}\label{sec:intro}

\IEEEPARstart{A}~high-performance controller is essential to the success of integrating large-scale wind energy into future power systems. Extensive investigations and recent lessons learned have confirmed that the variable and intermittent nature of wind indeed poses serious threats to both the reliability of power systems and the economic viability of wind energy \cite{UCTE06}. To minimize these threats, a wind turbine controller should not only maximize the amount of active power captured in a so-called {\em Maximum Power Tracking} (MPT) mode in normal situations, it should also allow the power captured be continuously regulated at a desired level in a so-called {\em Power Regulation} (PR) mode (other than clipping the power based on received instructions) when there is a system contingency. In addition, the controller should enable seamless switching between the MPT and PR modes, as well as maintain a desired power factor by also controlling the reactive power output.

One of the challenges facing the development of such a high-performance controller is the fact that the aerodynamic and mechanical parameters of a wind turbine are inherently uncertain, due to modeling and measurement errors, unknown optimal operating points, and other, possibly time-varying, ambient factors. For example, the $C_p$-surface of a wind turbine, which characterizes the amount of mechanical energy converted from the wind, is typically assumed to be known---or, at least, its optimal points are assumed to be known---in many existing controller designs. Unfortunately, such a surface is an empirical, statistical approximation, obtained based on up to three months of continuous experiment \cite{Riso01}. Thus, the $C_p$-surface may not be precisely known for control purposes. Other factors, such as changes in air density due to weather, variations in friction under different operating conditions, and measurement errors due to anemometer location, also contribute to the uncertainties. Indeed, a report from the National Renewable Energy Laboratory (NREL) \cite{NREL02}, which describes the result of its long-term test on different controllers operating on real, MW-level wind turbines, shows that the impact of these uncertainties on controller performance is significant and should be accounted for in controller design.

Another major challenge facing the development of such a high-performance controller is the fact that the mechanical and electrical parts of a modern wind turbine, which uses a Doubly Fed Induction Generator (DFIG), are tightly coupled. As will be detailed below, most studies have adopted a standard approach in the analysis and control of synchronous electric machines, in which the active and reactive powers are considered decoupled. With this approach, the active and reactive powers are adjusted via control of the mechanical and electrical parts, respectively, independent of each other. However, although a DFIG has some features of a synchronous machine, it is by nature an induction machine that exhibits strong electromechanical coupling among its rotor excitation current, rotor angular velocity, and electromagnetic torque. Hence, for performance reasons, the mechanical and electrical parts of a wind turbine with a DFIG should be considered synergistically in controller design.

The current literature offers a large collection of wind turbine controllers, including \cite{Iyasere08,Johnson04b,Johnson06,Galdi08,Chedid00,Beltran08,GengH09,Muljadi01,Senjyu06,Stol01,Wright03,Tarnowski07,Ko07,Marinescu04,ZhiDW07,WuF08,Pena96,Hopfensperger00,Monroy08,Peresada04,Hansen04}. However, as was alluded to above, most of the existing publications considered the mechanical and electrical parts separately (e.g., \cite{Iyasere08,Chedid00,Johnson04b,Johnson06,Beltran08,Galdi08,GengH09,Muljadi01,Senjyu06,Stol01,Wright03} considered only the former, while \cite{Ko07,ZhiDW07,Marinescu04,Tarnowski07,WuF08,Pena96,Monroy08,Peresada04,Hopfensperger00} considered only the latter), and for a few of those (e.g., \cite{Hansen04}) that considered both parts, its controller can only operate in the MPT mode, maximizing wind energy conversion, as opposed to both the MPT and PR modes. Moreover, although the existing work has provided valuable understanding in the control of wind turbines, only a few publications have addressed the issue of uncertainties. For example, \cite{Iyasere08,Johnson04b,Johnson06} proposed adaptive frameworks for controlling the mechanical part of wind turbines, so that the power captured is maximized, despite not knowing the $C_p$-surface.

In our recent work \cite{Tang09}, we developed a nonlinear controller that simultaneously enables control of the active power in both the MPT and PR modes, seamless switching between the two, and control of the reactive power so that a desirable power factor is maintained. These objectives were achieved by adjusting the rotor voltages of the electrical part and the blade pitch angle of the mechanical part, where the coupling between the two parts were taken into account in the controller design. Like most of the existing work, however, the controller in \cite{Tang09} assumed that the aerodynamic and mechanical parameters were known.

In this paper, we develop a controller that achieves such objectives and, at the same time, addresses the two aforementioned challenges, on uncertainties in the aerodynamic and mechanical parameters, and coupling between the mechanical and electrical parts. For the former, we show that the parametric uncertainties can be lumped into a scalar term, estimated via an uncertainty estimator in an inner loop, and circumvented in an outer, gradient-based minimization loop. For the latter, we show that the electromechanical coupling can be eliminated via feedback linearization on the electrical dynamics, following ideas from \cite{Tang09}. Finally, we demonstrate the effectiveness of the controller developed through simulation with a realistic wind profile from a wind farm in Oklahoma.

The outline of this paper is as follows: Section~\ref{sec:mod} models the wind turbine and formulates the problem. Section~\ref{sec:contro} describes the proposed controller. Section~\ref{sec:simres} presents the simulation results. Finally, Section~\ref{sec:concl} concludes the paper. The proof of the main theorem is included in the Appendix.

\section{Modeling and Problem Formulation}\label{sec:mod}

Consider a variable-speed wind turbine with a Doubly Fed Induction Generator (DFIG). The wind turbine consists of an electrical part and a mechanical part, the dynamics of which may be modeled as follows:

The dynamics of the electrical part in the $dq$ frame are described by a fourth-order state space model \cite{Bose02,Fadaeinedjad08}
\begin{align}
\begin{bmatrix}
\dot{\varphi}_{ds} \\
\dot{\varphi}_{qs} \\
\dot{\varphi}_{dr} \\
\dot{\varphi}_{qr}
\end{bmatrix}
&=\underbrace{
\begin{bmatrix}
-\frac{R_s}{\sigma L_s} & \omega_s & \frac{R_sL_m}{\sigma L_sL_r} & 0 \\
-\omega_s & -\frac{R_s}{\sigma L_s} & 0 & \frac{R_sL_m}{\sigma L_sL_r}\\
\frac{R_rL_m}{\sigma L_sL_r} & 0 & -\frac{R_r}{\sigma L_r} & \omega_s \\
0 & \frac{R_rL_m}{\sigma L_sL_r} & -\omega_s & -\frac{R_r}{\sigma L_r}
\end{bmatrix}}_{A}
\begin{bmatrix}
\varphi_{ds} \\
\varphi_{qs} \\
\varphi_{dr} \\
\varphi_{qr}
\end{bmatrix} \nonumber\displaybreak[0]\\
&+\underbrace{ \begin{bmatrix}
0 & 0 \\
0 & 0 \\
1 & 0 \\
0 & 1
\end{bmatrix}}_{B}
\begin{bmatrix}
v_{dr} \\
v_{qr}
\end{bmatrix} +
\begin{bmatrix}
v_{ds} \\
v_{qs} \\
-\omega_r\varphi_{qr} \\
\omega_r\varphi_{dr}
\end{bmatrix}, \label{eqn:electricaldynamics}
\end{align}
where $\varphi_{ds}$, $\varphi_{qs}$, $\varphi_{dr}$, $\varphi_{qr}\in \mathbb{R}$ are state variables representing the stator and rotor fluxes, $v_{dr}$, $v_{qr}\in \mathbb{R}$ are control variables representing the rotor voltages, $v_{ds}$, $v_{qs}\in \mathbb{R}$ are the constant stator voltages (that are not simultaneously zero), $\omega_s>0$ is the constant angular velocity of the synchronously rotating reference frame, $\omega_r>0$ is the rotor angular velocity, $R_s$, $R_r$ are the stator and rotor resistances, $L_s$, $L_r$, $L_m$ are the stator, rotor, and mutual inductances satisfying $L_s>L_m$ and $L_r>L_m$, $\sigma=1-\frac{L_m^2}{L_sL_r}$ is the leakage coefficient, and $A$, $B$ are constant matrices. In addition, the fluxes can be written as \cite{Fadaeinedjad08}
\begin{align}
\underbrace{
\begin{bmatrix}
\varphi_{ds} \\
\varphi_{qs} \\
\varphi_{dr} \\
\varphi_{qr}
\end{bmatrix}}_{\varphi}=
\begin{bmatrix}
L_s & 0 & L_m & 0\\
0 & L_s & 0 & L_m\\
L_m & 0 & L_r & 0\\
0 & L_m & 0 & L_r
\end{bmatrix} \underbrace{\begin{bmatrix}
i_{ds} \\
i_{qs} \\
i_{dr} \\
i_{qr}
\end{bmatrix}}_{i}, \label{eqn:phi}
\end{align}
where $i_{ds}$, $i_{qs}$, $i_{dr}$, $i_{qr}\in \mathbb{R}$ are the stator and rotor currents, while $\varphi=[\varphi_{ds}\;\varphi_{qs}\;\varphi_{dr}\;\varphi_{qr}]^T$ and $i=[i_{ds}\; i_{qs}\; i_{dr}\; i_{qr}]^T$ are introduced just for convenience. Furthermore, the active and reactive stator and rotor powers are given by \cite{LeiYZ06}
\begin{align}
P_s&=-v_{ds}i_{ds}-v_{qs}i_{qs},\quad Q_s=-v_{qs}i_{ds}+v_{ds}i_{qs},\label{eqn:Ps}\\
P_r&=-v_{dr}i_{dr}-v_{qr}i_{qr},\quad Q_r=-v_{qr}i_{dr}+v_{dr}i_{qr},\label{eqn:Pr}
\end{align}
and the total active and reactive powers of the turbine are
\begin{align}
P=P_s+P_r,\quad Q=Q_s+Q_r,\label{eqn:Q}
\end{align}
where positive (negative) values of $P$ and $Q$ mean that the turbine injects power into (draws power from) the grid.

The dynamics of the mechanical part are described by a first-order state space model \cite{Bose02}
\begin{align}
J\dot{\omega}_r=T_m-T_e-C_f\omega_r,\label{eqn:mech}
\end{align}
where the rotor angular velocity $\omega_r$ is another state variable, $J$ is the moment of inertia, $C_f$ is the friction coefficient, $T_m$ is the mechanical torque, and $T_e$ is the electromagnetic torque given by \cite{LeiYZ06}
\begin{align}
T_e=\varphi_{qs}i_{ds}-\varphi_{ds}i_{qs},\label{eqn:Te}
\end{align}
where positive (negative) value of $T_e$ means that the turbine acts as a generator (motor). The mechanical power captured by the wind turbine is \cite{Bianchi07}
\begin{align}
P_m=T_m\omega_r=\frac{1}{2}\rho AC_p(\lambda,\beta)V_{w}^3,\label{eqn:Pm}
\end{align}
where $\rho$ is the air density, $A=\pi R^2$ is the area swept by the rotor blades of radius $R$, $V_{w}$ is the wind speed, and $C_p(\lambda,\beta)$, commonly referred to as the $C_p$-surface, is the performance coefficient of the wind turbine, whose value is a function of the tip speed ratio $\lambda\in(0,\infty)$, defined as
\begin{align}
\lambda =\frac{\omega_rR}{V_{w}},\label{eqn:lambda}
\end{align}
and the blade pitch angle $\beta\in [\beta_{\min},\beta_{\max}]$, which is another control variable.

In order for results of this paper to be applicable to a broad class of wind turbines, no specific expression of $C_p(\lambda, \beta)$ will be assumed. Instead, $C_p(\lambda,\beta)$ will only be assumed to satisfy the following mild conditions for the purpose of analysis:
\begin{enumerate}
\renewcommand{\theenumi}{(A\arabic{enumi})}\itemsep-\parsep
\renewcommand{\labelenumi}{\theenumi}
\item Function $C_p(\lambda,\beta)$ is continuously differentiable in both $\lambda$ and $\beta$ over $\lambda\in(0,\infty)$ and $\beta\in[\beta_{\min}, \beta_{\max}]$.
\item There exists $c\in(0,\infty)$ such that for all $\lambda\in(0,\infty)$ and $\beta\in[\beta_{\min}, \beta_{\max}]$, we have $C_p(\lambda,\beta)\leq c\lambda$. This condition is mild because it is equivalent to saying that the mechanical torque $T_m$ is bounded from above, since $T_m\propto\frac{C_p(\lambda,\beta)}{\lambda}$ according to \eqref{eqn:Pm} and \eqref{eqn:lambda}.
\item For each fixed $\beta\in[\beta_{\min}, \beta_{\max}]$, there exists $\lambda_1\in(0, \infty)$ such that for all $\lambda\in(0, \lambda_1)$, we have $C_p(\lambda, \beta)>0$. This condition is also mild because turbines are designed to capture wind power over a wide range of $\lambda$, including times when $\lambda$ is small.
\item There exist $\underline{c}\in(-\infty,0)$ and $\overline{c}\in(0,\infty)$ such that for all $\lambda\in(0, \infty)$ and $\beta\in[\beta_{\min}, \beta_{\max}]$, we have $\underline{c}\leq\frac{\partial}{\partial \lambda}(\frac{C_p(\lambda, \beta)}{\lambda})\leq\overline{c}$.
\end{enumerate}

\begin{figure*}[tb]
\centering\includegraphics[width=0.9\textwidth]{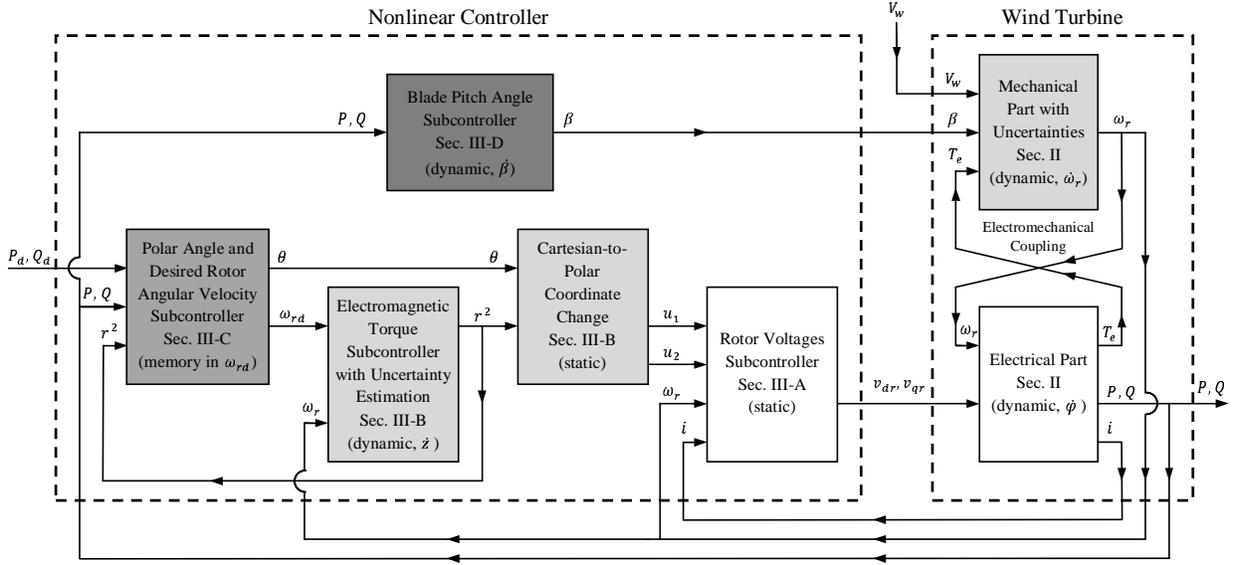}
\caption{Model of the wind turbine and architecture of the nonlinear controller.}
\label{fig:blockdiagram}
\end{figure*}

As it follows from the above, the wind turbine is modeled as a fifth-order, electromechanically-coupled, nonlinear system with state equations \eqref{eqn:electricaldynamics} and \eqref{eqn:mech}, output equations \eqref{eqn:Ps}--\eqref{eqn:Q}, state variables $\varphi_{ds}$, $\varphi_{qs}$, $\varphi_{dr}$, $\varphi_{qr}$, and $\omega_r$, control variables $v_{dr}$, $v_{qr}$, and $\beta$, output variables $P$ and $Q$, and exogenous ``disturbance'' $V_w$. A block diagram of this system is shown on the right-hand side of Figure~\ref{fig:blockdiagram}, in which the electromechanical coupling can be seen.

Given the above model, the problem addressed in this paper is: design a feedback controller, so that the active and reactive powers $P$ and $Q$ closely track some desired, possibly time-varying references $P_d$ and $Q_d$, assumed to be provided by a wind farm operator. When $P_d$ is larger than what the wind turbine is capable of generating, it means that the operator wants the turbine to operate in the {\em Maximum Power Tracking} (MPT) mode; otherwise, the {\em Power Regulation} (PR) mode is sought. By also providing $Q_d$, the operator indirectly specifies a desired power factor $\text{PF}_d=\frac{P_d}{\sqrt{P_d^2+Q_d^2}}$, around which the actual power factor $\text{PF}=\frac{P}{\sqrt{P^2+Q^2}}$ should be regulated. The controller may use $i$, $\omega_r$, $P$, and $Q$, which are all measurable, as feedback. The fluxes $\varphi$ may also be viewed as feedback, since they are bijectively related to $i$ through \eqref{eqn:phi}. Moreover, the controller may use values of all the electrical parameters (i.e., $\omega_s$, $R_s$, $R_r$, $L_s$, $L_r$, $L_m$, $v_{ds}$, and $v_{qs}$) and turbine-geometry-dependent parameters (i.e., $J$, $A$, $R$, $\beta_{\min}$, and $\beta_{\max}$), since these values are typically quite accurately known. However, it may not use values of the $C_p$-surface, the air density $\rho$, and the friction coefficient $C_f$, since these values are inherently uncertain and can change over time. Furthermore, the controller should not rely on the wind speed $V_w$, since it may not be accurately measured.

\section{Controller Design}\label{sec:contro}

In this section, we address the aforementioned problem by developing a nonlinear controller consisting of four subcontrollers. Figure~\ref{fig:blockdiagram} shows the architecture of the nonlinear controller, where each block represents a subcontroller. Note that the controller accepts $P_d$ and $Q_d$ as reference inputs, uses $i$, $\omega_r$, $P$, and $Q$ as feedback, and produces $v_{dr}$, $v_{qr}$, and $\beta$ as control inputs to the wind turbine. Moreover, the different gray levels of the blocks in Figure~\ref{fig:blockdiagram} represent our intended time-scale separation in the closed-loop dynamics: the darker a block, the slower its dynamics. The subcontrollers will be described in Sections~\ref{ssec:rotvolsubcon}--\ref{ssec:bladesubcon}. Note that Sections \ref{ssec:rotvolsubcon} and \ref{ssec:torquesubcon} up to the coordinate change are similar to our previous work \cite{Tang09}, while the rest of the paper contains new, unpublished results.

\subsection{Rotor Voltages Subcontroller}\label{ssec:rotvolsubcon}

Observe that although the electrical dynamics \eqref{eqn:electricaldynamics} are nonlinear, they possess a nice structure: the first and second rows of \eqref{eqn:electricaldynamics} are affine, consisting of linear terms and the constants $v_{ds}$ and $v_{qs}$, while the third and fourth are nonlinear, consisting of linear terms, the control variables $v_{dr}$ and $v_{qr}$, and the nonlinearities $-\omega_r\varphi_{qr}$ and $\omega_r\varphi_{dr}$ induced by the electromechanical coupling. Since the nonlinearities enter the dynamics the same way the control variables $v_{dr}$ and $v_{qr}$ do, we may use {\em feedback linearization} \cite{Khalil01} to cancel them and perform {\em pole placement} \cite{Chen99}, i.e., let
\begin{align}
v_{dr}&=\omega_r\varphi_{qr}-K_1^T\varphi+u_1, \label{eqn:vdr1}\\
v_{qr}&=-\omega_r\varphi_{dr}-K_2^T\varphi+u_2,\label{eqn:vqr1}
\end{align}
where $\omega_r\varphi_{qr}$ and $-\omega_r\varphi_{dr}$ are intended to cancel the nonlinearities, $-K_1^T\varphi$ and $-K_2^T\varphi$ with $K_1, K_2 \in\mathbb{R}^4$ are for pole placement, and $u_1$ and $u_2$ are new control variables to be designed in Section~\ref{ssec:torquesubcon}.

Substituting \eqref{eqn:vdr1} and \eqref{eqn:vqr1} into \eqref{eqn:electricaldynamics}, we get
\begin{align}
\dot{\varphi}=(A-BK)\varphi+\begin{bmatrix}
v_{ds}\;v_{qs}\;u_1\;u_2
\end{bmatrix}^T, \label{eqn:dotx}
\end{align}
where $K=[K_1\quad K_2]^T$ is the state feedback gain matrix. Since the electrical dynamics are physically allowed to be much faster than the mechanicals, we may choose $K$ in \eqref{eqn:dotx} to be such that $A-BK$ is asymptotically stable with very fast eigenvalues. With $K$ chosen as such and with relatively slow-varying $u_1$ and $u_2$, the linear differential equation \eqref{eqn:dotx} may be approximated by a linear algebraic equation:
\begin{align}
\varphi=-(A-BK)^{-1}
\begin{bmatrix}
v_{ds}\;v_{qs}\;u_1\;u_2
\end{bmatrix}^T.\label{eqn:equilibrium}
\end{align}
Consequently, the fifth-order state equations \eqref{eqn:electricaldynamics} and \eqref{eqn:mech} may be approximated by the first-order state equation \eqref{eqn:mech} along with algebraic relationships \eqref{eqn:vdr1}, \eqref{eqn:vqr1}, and \eqref{eqn:equilibrium}. This approximation will be made in all subsequent development (but not in simulation).

Note that \eqref{eqn:phi}, \eqref{eqn:vdr1}, and \eqref{eqn:vqr1} describe the {\em Rotor Voltages Subcontroller} block in Figure~\ref{fig:blockdiagram}.

\subsection{Electromagnetic Torque Subcontroller with Uncertainty Estimation}\label{ssec:torquesubcon}

Having addressed the electrical dynamics, we now consider the mechanicals, where the goal is to construct a subcontroller, which makes the rotor angular velocity $\omega_r$ track a desired, slow-varying reference $\omega_{rd}$, despite not knowing the aerodynamic and mechanical parameters listed at the end of Section~\ref{sec:mod}.

To come up with such a subcontroller, we first introduce a coordinate change. As was shown in our previous work \cite{Tang09}, because of \eqref{eqn:phi}, \eqref{eqn:Te}, and \eqref{eqn:equilibrium}, the electromagnetic torque $T_e$ may be expressed as a quadratic function of the new control variables $u_1$ and $u_2$, i.e.,
\begin{align}
T_e&=
\begin{bmatrix}
u_1 & u_2
\end{bmatrix}
\begin{bmatrix}
q_1 & q_2 \\ q_2 & q_3
\end{bmatrix}
\begin{bmatrix}
u_1 \\ u_2
\end{bmatrix}+
\begin{bmatrix}
b_1 & b_2
\end{bmatrix}
\begin{bmatrix}
u_1 \\ u_2
\end{bmatrix}
+a,\label{eqn:Tequa}
\end{align}
where $q_1$, $q_2$, $q_3$, $b_1$, $b_2$, and $a$ depend on the electrical parameters and the state feedback gain matrix $K$. Moreover, as was shown in \cite{Tang09}, this quadratic function is always {\em convex} because its associated Hessian matrix $\left[\begin{smallmatrix}
q_1 & q_2 \\
q_2 & q_3
\end{smallmatrix}\right]$
is always {\em positive definite}. Since the mechanical dynamics \eqref{eqn:mech}, in $\omega_r$, are driven by $T_e$, while $T_e$ in \eqref{eqn:Tequa} is a quadratic function of $u_1$ and $u_2$, the {\em two} new control variables $u_1$ and $u_2$ collectively affect {\em one} state variable $\omega_r$. This implies that there is a redundancy in $u_1$ and $u_2$, which may be exploited elsewhere. Since the quadratic function is always convex, this redundancy may be exposed via the following coordinate change \cite{Tang09}, which transforms $u_1, u_2 \in \mathbb{R}$ in a Cartesian coordinate system into $r \geq 0$ and $\theta \in [-\pi, \pi)$ in a polar coordinate system:
\begin{align}
r=\sqrt{z_1^2+z_2^2},\quad\theta=\operatorname{atan2}(z_2,z_1),\label{eqn:ztheta}
\end{align}
where
\begin{align}
\begin{bmatrix}
z_1 \\ z_2
\end{bmatrix}=D^{1/2}M^T\begin{bmatrix}
u_1 \\ u_2 \end{bmatrix}+\frac{1}{2}D^{-1/2}M^T\begin{bmatrix}
b_1 \\ b_2 \end{bmatrix}, \label{eqn:zfunu}
\end{align}
$\operatorname{atan2}()$ denotes the four-quadrant arctangent function, and $M$ and $D$ contain the eigenvectors and eigenvalues of $\left[\begin{smallmatrix}
q_1 & q_2 \\
q_2 & q_3
\end{smallmatrix}\right]$ on their columns and diagonal, respectively, i.e., $M^T\left[\begin{smallmatrix}
q_1 & q_2 \\
q_2 & q_3
\end{smallmatrix}\right]M=D$. In the polar coordinates, it follows from \eqref{eqn:Tequa}--\eqref{eqn:zfunu} that
\begin{align}
T_e&=r^2+a',\label{eqn:Techange2}
\end{align}
where $a'=-\frac{v_{ds}^2+v_{qs}^2}{4\omega_sR_s}$ is always negative. From \eqref{eqn:mech} and \eqref{eqn:Techange2}, we see that in the polar coordinates, $r^2$ is responsible for driving the mechanical dynamics in $\omega_r$ and, hence, may be viewed as an {\em equivalent electromagnetic torque}, differed from $T_e$ only by a constant $a'$. On the other hand, the {\em polar angle} $\theta$ has no impact on the mechanical dynamics and, thus, represents the redundancy that will be exploited later, in Section~\ref{ssec:anglevelsubcon}.

Note that \eqref{eqn:ztheta} and \eqref{eqn:zfunu} describe the {\em Cartesian-to-Polar Coordinate Change} block in Figure~\ref{fig:blockdiagram}.

Having introduced the coordinate change, we next show that the unknown aerodynamic and mechanical parameters, listed at the end of Section~\ref{sec:mod}, can be lumped into a scalar term, simplifying the problem. Combining \eqref{eqn:mech}, \eqref{eqn:Pm}, \eqref{eqn:lambda}, and \eqref{eqn:Techange2},
\begin{align}
J\dot{\omega}_r=\frac{\frac{1}{2}\rho AC_p(\frac{\omega_rR}{V_w}, \beta)V_w^3}{\omega_r}-r^2-a'-C_f\omega_r. \label{eqn:mechTechange}
\end{align}
Notice that the unknown parameters---namely, the $C_p$-surface, the air density $\rho$, the friction coefficient $C_f$, and the wind speed $V_w$---all appear in \eqref{eqn:mechTechange}. Moreover, these unknown parameters can be separated from the ``control input'' $r^2$ and lumped into a scalar function $g(\omega_r, \beta, V_w)$, defined as
\begin{align}
g(\omega_r, \beta, V_w)=\frac{\frac{1}{2}\rho AC_p(\frac{\omega_rR}{V_w}, \beta)V_w^3}{\omega_r}-a'-C_f\omega_r.\label{eqn:g}
\end{align}
With $g(\omega_r, \beta, V_w)$ in \eqref{eqn:g} representing the aggregated uncertainties, the first-order dynamics \eqref{eqn:mechTechange} are simplified to
\begin{align}
\dot{\omega}_r=\frac{1}{J}(g(\omega_r, \beta, V_w)-r^2). \label{eqn:mechg}
\end{align}

To design a controller for $r^2$, which allows the rotor angular velocity $\omega_r$ to track a desired, slow-varying reference $\omega_{rd}$ despite the unknown scalar function $g(\omega_r, \beta, V_w)$, consider a first-order nonlinear system
\begin{align}
\dot{x}=\frac{1}{J}(f(x)+u),\label{eqn:simplemodel}
\end{align}
where $x\in\mathbb{R}$ is the state, $u\in\mathbb{R}$ is the input, and $f(x)$ is a {\em known function} of $x$. Obviously, to drive $x$ to some desired value $x_d\in\mathbb{R}$, we may apply feedback linearization \cite{Khalil01} to cancel $f(x)$ and insert linear dynamics, i.e., let
\begin{align}
u=-f(x)-\alpha(x-x_d),\label{eqn:controlu}
\end{align}
where $\alpha\in\mathbb{R}$ is the controller gain. Combining \eqref{eqn:simplemodel} with \eqref{eqn:controlu} yields the closed-loop dynamics
\begin{align}
\dot{x}=-\frac{\alpha}{J}(x-x_d).\label{eqn:desireddyn}
\end{align}
Thus, if $\alpha$ is positive, $x$ in \eqref{eqn:desireddyn} asymptotically goes to $x_d$.

Now suppose $f(x)$ in \eqref{eqn:simplemodel} is {\em unknown} but a {\em constant}, denoted simply as $f\in\mathbb{R}$ (we will relax the assumption that it is a constant shortly). With $f$ being unknown, the controller \eqref{eqn:controlu} is no longer applicable. To overcome this limitation, we may first introduce a reduced-order estimator \cite{Brogan91}, which calculates an estimate $\hat{f}\in\mathbb{R}$ of $f$, and then replace $f(x)$ in \eqref{eqn:controlu} by the estimate $\hat{f}$:
\begin{align}
\dot{z}&=-\frac{h}{J}(u+\hat{f}),\label{eqn:observerstate}\\
\hat{f}&=z+hx,\label{eqn:observeroutput}\\
u&=-\hat{f}-\alpha(x-x_d),\label{eqn:controlu2}
\end{align}
where $z\in\mathbb{R}$ is the estimator state and $h\in\mathbb{R}$ is the estimator gain. Defining the estimation error as $\tilde{f}=f-\hat{f}$ and combining \eqref{eqn:simplemodel} with \eqref{eqn:observerstate}--\eqref{eqn:controlu2} yield closed-loop dynamics
\begin{align}
\dot{\tilde{f}}&=-\dot{\hat{f}}=-\dot{z}-h\dot{x}=-\frac{h}{J}\tilde{f},\label{eqn:errdynf2}\\
\dot{x}&=\frac{1}{J}(f-\hat{f}-\alpha(x-x_d))=\frac{1}{J}(\tilde{f}-\alpha(x-x_d)).\label{eqn:desireddyn2}
\end{align}
Hence, by letting both $\alpha$ and $h$ be positive, both $\tilde{f}$ and $x$ in \eqref{eqn:errdynf2} and \eqref{eqn:desireddyn2} asymptotically go to $0$ and $x_d$, respectively.

Next, suppose both the state $x$ and the desired value $x_d$ must be {\em positive}, instead of being anywhere in $\mathbb{R}$. With this restriction, the controller with uncertainy estimation \eqref{eqn:observerstate}--\eqref{eqn:controlu2} needs to be modified, because for some initial conditions, it is possible that $x$ can become nonpositive. One way to modify the controller is to replace the linear term $x-x_d$ in \eqref{eqn:controlu2} by a logarithmic one $\ln\frac{x}{x_d}$, resulting in
\begin{align}
u=-\hat{f}-\alpha\ln\frac{x}{x_d}.\label{eqn:controlu3}
\end{align}
With \eqref{eqn:observerstate}, \eqref{eqn:observeroutput}, and \eqref{eqn:controlu3}, the closed-loop dynamics become
\begin{align}
\dot{\tilde{f}}&=-\frac{h}{J}\tilde{f},\label{eqn:errdynf3}\\
\dot{x}&=\frac{1}{J}(\tilde{f}-\alpha\ln\frac{x}{x_d}).\label{eqn:desireddyn3}
\end{align}
Note from \eqref{eqn:desireddyn3} that for any $\tilde{f}\in\mathbb{R}$, there exists positive $x$, sufficiently small, such that $\dot{x}$ is positive. Therefore, for any initial condition $(\tilde{f}(0), x(0))$ with positive $x(0)$, $x(t)$ will remain positive, suggesting that the modification \eqref{eqn:controlu3} satisfies the restriction.

Now suppose the input $u$ must be {\em nonpositive}. With this additional restriction, \eqref{eqn:controlu3} needs to be further modified. One way to do so is to force the right-hand side of \eqref{eqn:controlu3} to be nonpositive, leading to
\begin{align}
u=-\max\{\hat{f}+\alpha\ln\frac{x}{x_d},0\}.\label{eqn:controlu4}
\end{align}
Clearly, with \eqref{eqn:controlu4}, $u$ is always nonpositive.

Finally, suppose $f$ is an {\em unknown function} of $x$, denoted as $f(x)$. With this relaxation, we may associate the first-order nonlinear system \eqref{eqn:simplemodel} with the first-order dynamics \eqref{eqn:mechg} by viewing $x$ as $\omega_r$, $x_d$ as $\omega_{rd}$, $u$ as $-r^2$, $f(x)$ as $g(\omega_r, \beta, V_w)$ (treating $\beta$ and $V_w$ as constants), and $\hat{f}$ as $\hat{g}$ (i.e., $\hat{g}$ is an estimate of $g(\omega_r, \beta, V_w)$). Based on this association, \eqref{eqn:observerstate}, \eqref{eqn:observeroutput}, and \eqref{eqn:controlu4} can be written as
\begin{align}
\dot{z}&=-\frac{h}{J}(-r^2+\hat{g}),\label{eqn:newdyn}\\
\hat{g}&=z+h\omega_r,\label{eqn:ghat}\\
r^2&=\max\{\hat{g}+\alpha\ln\frac{\omega_r}{\omega_{rd}},0\}. \label{eq:new_r}
\end{align}

Having derived the controller with uncertainty estimation \eqref{eqn:newdyn}--\eqref{eq:new_r}, we now analyze its behavior. To do so, some setup is needed: first, suppose $\omega_{rd}$, $\beta$, and $V_w$ are constants. Second, as was shown in \cite{Tang09}, because of Assumptions (A1)--(A3) in Section~\ref{sec:mod}, there exists $\omega_r^{(1)}\in(0,\infty)$ such that $g(\omega_r^{(1)}, \beta, V_w)=0$ and $g(\omega_r, \beta, V_w)>0$ for all $\omega_r\in(0, \omega_r^{(1)})$. Third, using \eqref{eqn:lambda}, \eqref{eqn:g}, and Assumptions (A1) and (A4), it is straightforward to show that there exist $\underline{\gamma}\in(-\infty, 0)$ and $\overline{\gamma}\in(0, \infty)$ such that $\underline{\gamma}\leq\frac{\partial}{\partial \omega_r}g(\omega_r, \beta, V_w)\leq\overline{\gamma}$ for all $\omega_r\in(0, \infty)$. Finally, with \eqref{eqn:mechg} and \eqref{eqn:newdyn}--\eqref{eq:new_r} and with $(\omega_r, \hat{g})$ as state variables (instead of $(\omega_r, z)$), the closed-loop dynamics can be expressed as
\begin{align}
\dot{\omega}_r&=\frac{1}{J}(g(\omega_r, \beta, V_w)-\max\{\hat{g}+\alpha\ln\frac{\omega_r}{\omega_{rd}}, 0\}),\label{eqn:closedsys1}\\
\dot{\hat{g}}&=\dot{z}+h\dot{\omega}_r=\frac{h}{J}(g(\omega_r, \beta, V_w)-\hat{g}).\label{eqn:closedsys2}
\end{align}

The following theorem characterizes the stability properties of the closed-loop system \eqref{eqn:closedsys1} and \eqref{eqn:closedsys2}:

\begin{theorem}\label{thm:closedsys}{\em Consider the closed-loop system \eqref{eqn:closedsys1} and \eqref{eqn:closedsys2}. Suppose $\omega_{rd}$, $\beta$, and $V_w$ are constants with $0<\omega_{rd}\leq\omega_r^{(1)}$, where $\omega_r^{(1)}$, along with $\underline{\gamma}$ and $\overline{\gamma}$, is as defined above. Let $D=\{(\omega_r, \hat{g})|0<\omega_r\leq\omega_r^{(1)}, \hat{g}\in\mathbb{R}\}\subset\mathbb{R}^2$. If the controller gain $\alpha$ is positive and the estimator gain $h$ is sufficiently large, i.e.,
\begin{align}
\begin{array}{ll}
h>\overline{\gamma}&\text{if}\quad\overline{\gamma}\geq-\frac{1}{3}\underline{\gamma},\\
h>-\frac{(\overline{\gamma}-\underline{\gamma})^2}{8(\underline{\gamma}+\overline{\gamma})}&\text{otherwise},
\end{array}\label{eq:condh}
\end{align}
then: (i) the system has a unique equilibrium point at $(\omega_{rd}, g(\omega_{rd}, \beta, V_w))$ in $D$; (ii) the set $D$ is a positively invariant set, i.e., if $(\omega_r(0), \hat{g}(0))\in D$, then $(\omega_r(t), \hat{g}(t))\in D$ $\forall t\geq0$; and (iii) the equilibrium point $(\omega_{rd}, g(\omega_{rd}, \beta, V_w))$ is locally asymptotically stable with a domain of attraction $D$.}
\end{theorem}

\begin{IEEEproof}
See the Appendix.
\end{IEEEproof}

Theorem~\ref{thm:closedsys} says that, by using the electromagnetic torque subcontroller with uncertainty estimation \eqref{eqn:newdyn}--\eqref{eq:new_r}, if the gains $\alpha$ and $h$ are positive and sufficiently large and if the desired reference $\omega_{rd}$ does not exceed $\omega_r^{(1)}$, then the rotor angular velocity $\omega_r$ asymptotically converges to $\omega_{rd}$ if $\omega_{rd}$, $\beta$, and $V_w$ are constants and closely tracks $\omega_{rd}$ if they are slow-varying. Notice that the gains $\alpha$ and $h$ can be chosen independently of each other. Also, the condition ``$\omega_{rd}\leq\omega_r^{(1)}$'' is practically always satisfied, as $\omega_r^{(1)}$ is extremely large \cite{Tang09}.

Note that \eqref{eqn:newdyn}--\eqref{eq:new_r} describe the {\em Electromagnetic Torque Subcontroller with Uncertainty Estimation} block in Figure~\ref{fig:blockdiagram}.

\subsection{Polar Angle and Desired Rotor Angular Velocity Subcontroller}\label{ssec:anglevelsubcon}

\begin{figure}[tb]
\centering\includegraphics[width=0.9\linewidth]{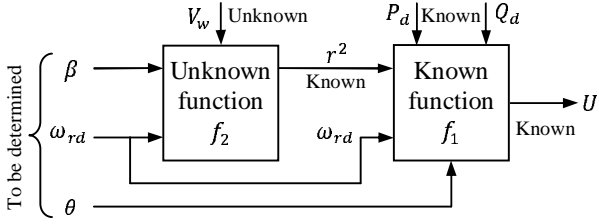}
\caption{Relationships among the performance measure $U$, the to-be-determined variables $\theta$, $\omega_{rd}$, and $\beta$, and the exogenous variables $V_w$, $P_d$, and $Q_d$.}
\label{fig:blocks}
\end{figure}

Up to this point in the paper, we have yet to specify how $\theta$, $\omega_{rd}$, and $\beta$ are determined. To do so, we first introduce a scalar performance measure and express this measure as a function of $\theta$, $\omega_{rd}$, and $\beta$. We then present a method for choosing these variables, which optimizes the measure.

Recall that the ultimate goal is to make the active and reactive powers $P$ and $Q$ track some desired references $P_d$ and $Q_d$ as closely as possible. Hence, it is useful to introduce a scalar performance measure, which characterizes how far $P$ and $Q$ are from $P_d$ and $Q_d$. One such measure, denoted as $U$, is given by
\begin{align}
U=\frac{1}{2}\begin{bmatrix} P-P_d & Q-Q_d \end{bmatrix}
\begin{bmatrix}
w_p & w_{pq} \\
w_{pq} & w_q
\end{bmatrix}
\begin{bmatrix}
P-P_d \\
Q-Q_d
\end{bmatrix},\label{eqn:objfun}
\end{align}
where $w_p$, $w_q$, and $w_{pq}$ are design parameters satisfying $w_p>0$ and $w_p w_q>w_{pq}^2$, so that
$\left[\begin{smallmatrix}
w_p & w_{pq} \\
w_{pq} & w_q
\end{smallmatrix}\right]$
is a positive definite matrix. With these design parameters, one may specify how the differences $P-P_d$ and $Q-Q_d$ and their product $(P-P_d)(Q-Q_d)$ are penalized. Moreover, with $U$ in \eqref{eqn:objfun} being a quadratic, positive definite function of $P-P_d$ and $Q-Q_d$, the smaller $U$ is, the better the ultimate goal is achieved.

Having defined the performance measure $U$ via \eqref{eqn:objfun}, we next establish the following statement: if the subcontrollers in Sections~\ref{ssec:rotvolsubcon} and~\ref{ssec:torquesubcon} are used with $K$ chosen so that $A-BK$ has very fast eigenvalues, $\alpha$ chosen to be positive, and $h$ chosen to satisfy \eqref{eq:condh}, and if $\theta$, $\omega_{rd}$, $\beta$, $V_w$, $P_{d}$, and $Q_{d}$ are all constants, then after a short transient, $U$ may be expressed as a {\em known} function $f_1$ of $r^2$, $\theta$, $\omega_{rd}$, $P_d$, and $Q_d$, while $r^2$, in turn, may be expressed as an {\em unknown} function $f_2$ of $\omega_{rd}$, $\beta$, and $V_w$, i.e.,
\begin{align}
U&=f_1(r^2, \theta, \omega_{rd}, P_d, Q_d), \label{eqn:Uf}\\
r^2&=f_2(\omega_{rd}, \beta, V_w), \label{eqn:f2}
\end{align}
as shown in Figure~\ref{fig:blocks}. To establish this statement, suppose the hypothesis is true. Then, after a short transient, it follows from \eqref{eqn:objfun} that $U$ is a known function of $P$, $Q$, $P_d$, and $Q_d$; from \eqref{eqn:phi}, \eqref{eqn:Ps}--\eqref{eqn:Q}, \eqref{eqn:vdr1}, and \eqref{eqn:vqr1} that $P$ and $Q$ are known functions of $\varphi$, $\omega_r$, $u_1$, and $u_2$; from \eqref{eqn:equilibrium} that $\varphi$ is a known function of $u_1$ and $u_2$; from \eqref{eqn:ztheta} and \eqref{eqn:zfunu} that $u_1$ and $u_2$ are known functions of $r^2$ and $\theta$; and from Theorem~\ref{thm:closedsys} that $\omega_r=\omega_{rd}$. Thus, \eqref{eqn:Uf} holds with $f_1$ being known. On the other hand, it follows from \eqref{eq:new_r} and Theorem~\ref{thm:closedsys} that $r^2=g(\omega_{rd}, \beta, V_w).$ Hence, \eqref{eqn:f2} holds with $f_2$ being unknown.

Equations \eqref{eqn:Uf} and \eqref{eqn:f2}, which are represented in Figure~\ref{fig:blocks}, suggest that $U$ is a function of the to-be-determined variables $\theta$, $\omega_{rd}$, and $\beta$ as well as the exogenous variables $V_w$, $P_d$, and $Q_d$. Given that the smaller $U$ is the better, these to-be-determined variables may be chosen to minimize $U$. However, such minimization is difficult to carry out because although $P_d$ and $Q_d$ are known, $V_w$ is not. To make matter worse, since $f_1$ is known but $f_2$ is not, the objective function is not entirely known. Somewhat fortunately, as was shown in Figure~\ref{fig:blocks}, $\theta$ affects $U$ only through $f_1$ and not $f_2$. Therefore, $\theta$ may be chosen to minimize $U$ for any given $r^2$, $\omega_{rd}$, $P_d$, and $Q_d$, i.e.,
\begin{align}
\theta=\operatorname{arg\,min}_{x\in[-\pi,\pi)}f_1(r^2,x,\omega_{rd}, P_d, Q_d), \label{eqn:thetacon}
\end{align}
which is implementable since $r^2$, $\omega_{rd}$, $P_d$, and $Q_d$ are all known. With $\theta$ chosen as in \eqref{eqn:thetacon}, the minimization problem reduces from a three-dimensional problem to a two-dimensional one, depending only on $\omega_{rd}$ and $\beta$. Since the objective function upon absorbing $\theta$ is unknown and since $V_w$ may change quickly, instead of minimizing $U$ with respect to both $\omega_{rd}$ and $\beta$---which may take a long time---we decide to sacrifice freedom for speed, minimizing $U$ only with respect to $\omega_{rd}$ and updating $\beta$ in a relatively slower fashion, which will be described in Section~\ref{ssec:bladesubcon}.

\begin{figure}[tb]
\centering\includegraphics[width=0.9\linewidth]{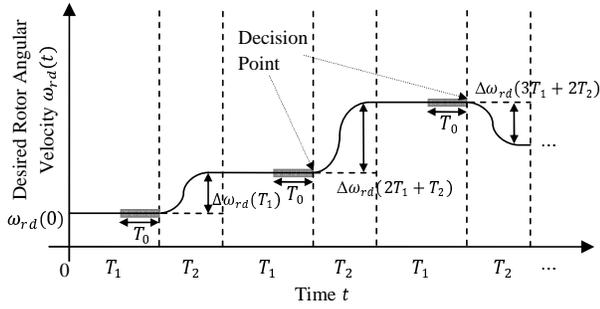}
\caption{A graphical illustration of the gradient-like approach.}
\label{fig:wrd}
\end{figure}

The minimization of $U$ with respect to $\omega_{rd}$ is carried out based on a gradient-like approach as shown in Figure~\ref{fig:wrd}. To explain the rationale behind this approach, suppose $\beta$, $V_w$, $P_d$, and $Q_d$ are constants. Then, according to \eqref{eqn:Uf}--\eqref{eqn:thetacon}, $U$ is an unknown function of $\omega_{rd}$. Because this function is not known, its gradient $\frac{\partial U}{\partial \omega_{rd}}$ at any $\omega_{rd}$ cannot be evaluated. To alleviate this issue, we evaluate $U$ at two nearby $\omega_{rd}$'s, use the two evaluated $U$'s to obtain an estimate of the gradient $\frac{\partial U}{\partial \omega_{rd}}$, and move $\omega_{rd}$ along the direction where $U$ decreases, by an amount which depends on the gradient estimate. This idea is illustrated in Figure~\ref{fig:wrd} and described precisely as follows: the desired rotor angular velocity $\omega_{rd}(t)$ is set to an initial value $\omega_{rd}(0)$ at time $t=0$ and held constant until $t=T_1$, where $T_1$ should be sufficiently large so that both the electrical and mechanical dynamics have a chance to reach steady-state, but not too large which causes the minimization to be too slow. From time $t=T_1-T_0$ to $t=T_1$, the average of $U(t)$, i.e., $\frac{1}{T_0}\int_{T_1-T_0}^{T_1}U(t)dt$, is recorded as the first value needed to obtain a gradient estimate. Similar to $T_1$, $T_0$ should be large enough so that small fluctuations in $U(t)$ (induced perhaps by a noisy $V_w$) are averaged out, but not too large which causes transient in the dynamics to be included. The variable $\omega_{rd}(t)$ is then changed gradually in an S-shape manner from $\omega_{rd}(0)$ at time $t=T_1$ to a nearby $\omega_{rd}(0)+\Delta\omega_{rd}(T_1)$ at $t=T_1+T_2$, where $\Delta\omega_{rd}(T_1)$ is an initial stepsize, and $T_2$ should be sufficiently large but not overly so, so that the transition in $\omega_{rd}(t)$ is smooth and yet not too slow. The variable $\omega_{rd}(t)$ is then held constant until $t=2T_1+T_2$, and the average of $U(t)$ from $t=2T_1+T_2-T_0$ to $t=2T_1+T_2$, i.e., $\frac{1}{T_0}\int_{2T_1+T_2-T_0}^{2T_1+T_2}U(t)dt$, is recorded as the second value needed to obtain the gradient estimate. At time $t=2T_1+T_2$, the two recorded values are used to form the gradient estimate, which is in turn used to decide a new stepsize $\Delta\omega_{rd}(2T_1+T_2)$ through
\begin{align}
&\Delta\omega_{rd}(2T_1+T_2)=\nonumber\\
&\;-\epsilon_1\sat\left(\frac{\frac{1}{T_0}\int_{2T_1+T_2-T_0}^{2T_1+T_2}U(t)dt-\frac{1}{T_0}\int_{T_1-T_0}^{T_1}U(t)dt}{\epsilon_2\Delta\omega_{rd}(T_1)}\right), \label{eqn:deltawrd}
\end{align}
where $\epsilon_1>0$ and $\epsilon_2>0$ are design parameters that define the new stepsize $\Delta\omega_{rd}(2T_1+T_2)$, and $\operatorname{sat}()$ denotes the standard saturation function that limits $\Delta\omega_{rd}(2T_1+T_2)$ to $\pm\epsilon_1$. Upon deciding $\Delta\omega_{rd}(2T_1+T_2)$, $\omega_{rd}(t)$ is again changed in an S-shape manner from $\omega_{rd}(0)+\Delta\omega_{rd}(T_1)$ at $t=2T_1+T_2$ to $\omega_{rd}(0)+\Delta\omega_{rd}(T_1)+\Delta\omega_{rd}(2T_1+T_2)$ at $t=2T_1+2T_2$, in a way similar to the time interval $[T_1,T_1+T_2]$. The process then repeats with the second recorded value from the previous cycle $[0,2T_1+T_2]$ becoming the first recorded value for the next cycle $[T_1+T_2,3T_1+2T_2]$, and so on. Therefore, with this gradient-like approach, $\omega_{rd}$ is guaranteed to approach a local minimum when $\beta$, $V_w$, $P_d$, and $Q_d$ are constants, and track a local minimum when they are slow-varying.

Note that \eqref{eqn:objfun}, \eqref{eqn:thetacon}, and \eqref{eqn:deltawrd} describe the {\em Polar Angle and Desired Rotor Angular Velocity Subcontroller} block in Figure~\ref{fig:blockdiagram}.

\subsection{Blade Pitch Angle Subcontroller}\label{ssec:bladesubcon}

As was mentioned, in order to speed up the minimization, we have decided to minimize $U$ only with respect to $\omega_{rd}$, leaving the blade pitch angle $\beta$ as the remaining undetermined variable. Given that an active power $P$ that is larger than the rated value $P_{\text{rated}}$ of the turbine may cause damage, we decide to use $\beta$ to prevent $P$ from exceeding $P_{\text{rated}}$, thereby protecting the turbine. Specifically, we let $\beta$ be updated according to
\begin{align}
\dot{\beta}=
\begin{cases}
0 & \text{if}\ \beta=\beta_{\min}\ \text{and}\ P<P_{\text{rated}},\\
0 & \text{if}\ \beta=\beta_{\max}\ \text{and}\ P>P_{\text{rated}},\\
-{\epsilon_3} (P_{\text{rated}}-P) & \text{otherwise},
\end{cases}
\label{eqn:betacont}
\end{align}
where $\epsilon_3>0$ is a design parameter that dictates the rate at which $\beta$ changes. Note that with \eqref{eqn:betacont}, $\beta$ is guaranteed to lie between $\beta_{\min}$ and $\beta_{\max}$. Moreover, when $P$ is above (below) $P_{\text{rated}}$, $\beta$ increases (decreases) if possible, in order to try to capture less (more) wind power, which leads to a smaller (larger) $P$.

Note that \eqref{eqn:betacont} describes the {\em Blade Pitch Angle Subcontroller} block in Figure~\ref{fig:blockdiagram}.

\begin{remark}The blade pitch angle subcontroller may be designed based on other considerations. For example, if the forecast of, say, the hourly-average wind speed $\overline{V}_w$ is available, for blade protection $\beta$ may be chosen as $\beta=F(\overline{V}_w)$ for some non-decreasing function $F:(0,\infty)\rightarrow{[\beta_{\min}, \beta_{\max}]}$.
\end{remark}

\section{Simulation Results}\label{sec:simres}

To demonstrate the capability and effectiveness of the proposed controller, simulation has been carried out in MATLAB. To describe the simulation settings and results, both the per-unit and physical unit systems will be used interchangeably.

The simulation settings are as follows: we consider a $1.5\operatorname{MW}$, $575\operatorname{V}$, $60\operatorname{Hz}$ wind turbine that is essentially adopted from the Distributed Resources Library in MATLAB/Simulink R2007a. The values of the wind turbine parameters are: $\omega_s=1 \operatorname{pu}$, $R_s=0.00706 \operatorname{pu}$, $R_r=0.005 \operatorname{pu}$, $L_s=3.071 \operatorname{pu}$, $L_r=3.056 \operatorname{pu}$, $L_m=2.9 \operatorname{pu}$, $v_{ds}=1 \operatorname{pu}$, $v_{qs}=0 \operatorname{pu}$, $J=10.08 \operatorname{pu}$, $A=4656.6 \operatorname{m^2}$, $R=38.5 \operatorname{m}$, $\beta_{\min}=0\operatorname{deg}$, $\beta_{\max}=30\operatorname{deg}$, and $C_f=0.01 \operatorname{pu}$. The $C_p$-surface adopted by MATLAB, which is taken from \cite{Heier98}, is $C_p(\lambda,\beta)=c_1\Bigl(\frac{c_2}{\lambda_i}-c_3\beta-c_4\Bigl)e^{\frac{-c_5}{\lambda_i}}+c_6\lambda$, where $\frac{1}{\lambda_i}=\frac{1}{\lambda+0.08\beta}-\frac{0.035}{\beta^3+1}$, $c_1=0.5176$, $c_2=116$, $c_3=0.4$, $c_4=5$, $c_5=21$, and $c_6=0.0068$. The mechanical power captured by the wind turbine is $P_m(\operatorname{pu})=\frac{P_\text{nom}P_{\text{wind}\_\text{base}}}{P_{\text{elec}\_\text{base}}}\,C_p(\operatorname{pu})\,V_w(\operatorname{pu})^3$, where $P_m(\operatorname{pu})=\frac{P_m}{P_{\text{nom}}}$, $P_{\text{nom}}=1.5\operatorname{MW}$ is the nominal mechanical power, $P_{\text{wind}\_\text{base}}=0.73\operatorname{pu}$ is the maximum power at the base wind speed, $P_{\text{elec}\_\text{base}}=1.5\times10^6/0.9\operatorname{VA}$ is the base power of the electrical generator, $C_p(\operatorname{pu})=\frac{C_p}{C_{p\_\text{nom}}}$, $C_{p\_\text{nom}}=0.48$ is the peak of the $C_p$-surface, $V_w(\operatorname{pu})=\frac{V_w}{V_{w\_\text{base}}}$, and $V_{w\_\text{base}}=12\operatorname{m/s}$ is the base wind speed. Note that the maximum mechanical power, captured at the base wind speed, is $0.657 \operatorname{pu}$. The tip speed ratio is $\lambda(\operatorname{pu})=\frac{\frac{\omega_r(\operatorname{pu})}{\omega_{r\_\text{base}}}}{V_w(\operatorname{pu})}$, where $\lambda(\operatorname{pu})=\frac{\lambda}{\lambda_{\text{nom}}}$, $\lambda_{\text{nom}}=8.1$ is the $\lambda$ that yields the peak of the $C_p$-surface, $\omega_{r\_\text{base}}=1.2\operatorname{pu}$ is the base rotational speed, $\omega_r(\operatorname{pu})=\frac{\omega_r}{\omega_{r\_\text{nom}}}$, and $\omega_{r\_\text{nom}}=2.1039\operatorname{rad/sec}$ is the nominal rotor angular velocity. For more details on these parameters and values, see the MATLAB documentation.

As for the proposed controller, we choose its parameters as follows: for the Rotor Voltages Subcontroller, we let the desired closed-loop eigenvalues of the electrical dynamics be at $-10$, $-15$, and $-20\pm5j$. Using MATLAB's {\tt place()} function, the state feedback gain matrix $K=[K_1\;K_2]^T$ that yields these eigenvalues is found to be $K_1=[12277\;\mbox{-4493.8}\;32.7\;\mbox{-6.3}]^T$ and $K_2=[\mbox{-1615.4}\;12117\;\mbox{-0.4}\;32.2]^T$. Moreover, we let $\alpha=5$ and $h=16$ for the Electromagnetic Torque Subcontroller with Uncertainty Estimation; $w_p=10$, $w_q=1$, $w_{pq}=0$, $\epsilon_1=0.025$, $\epsilon_2=2$, $T_0=1\operatorname{s}$, $T_1=4\operatorname{s}$, and $T_2=6\operatorname{s}$ for the Polar Angle and Desired Rotor Angular Velocity Subcontroller; and $\epsilon_3=2.7$ and $P_{\text{rated}}=1\operatorname{pu}$ for the Blade Pitch Angle Subcontroller.

\begin{figure*}[tb]
\centering\includegraphics[width=0.95\textwidth]{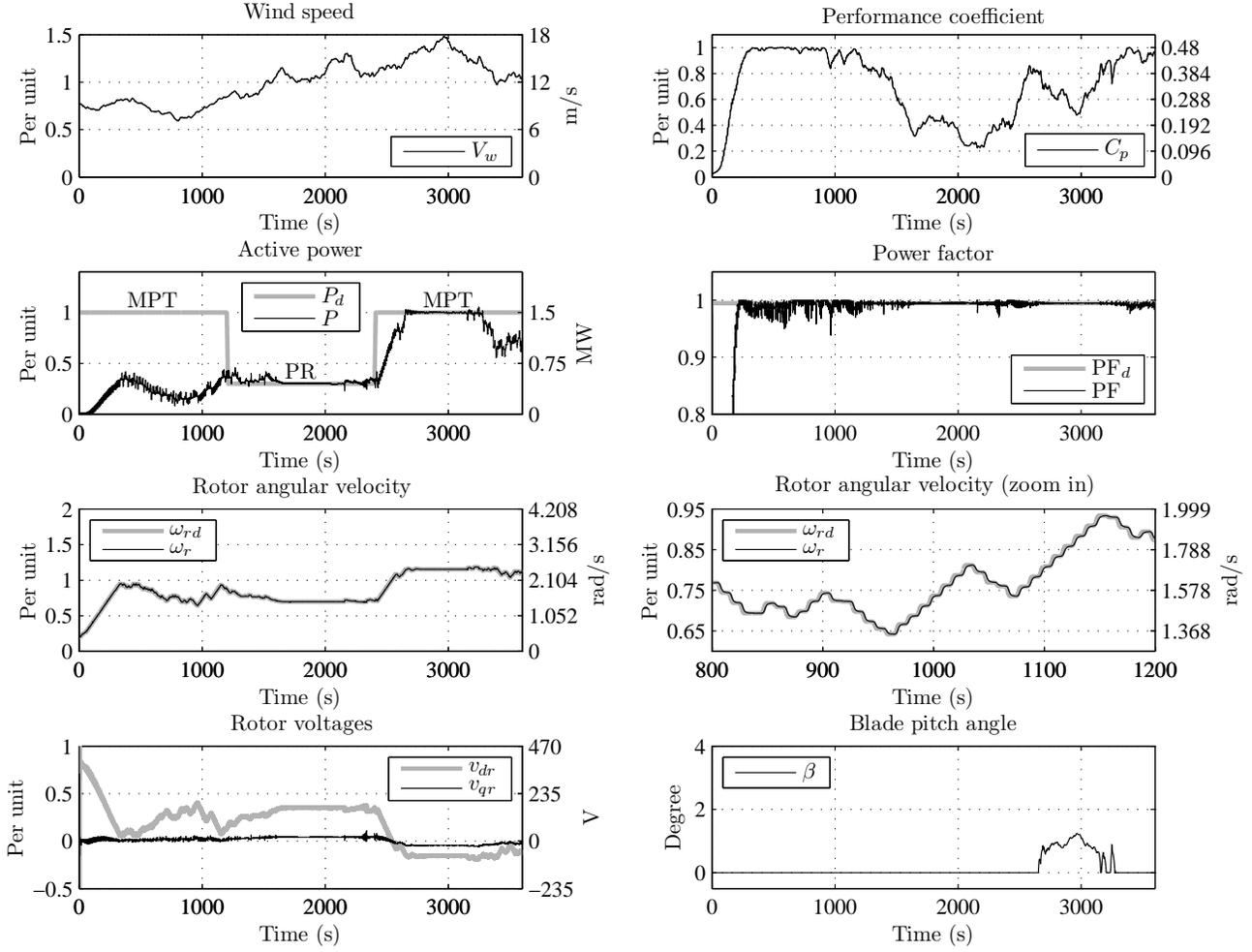}
\caption{Effective operation in both the MPT and PR modes and seamless switching between them under an actual wind profile from a wind farm located in northwest Oklahoma.}
\label{fig:sim11}
\end{figure*}

The simulation results are as follows: we consider a scenario where the wind speed $V_w$ is derived from an actual wind profile from a wind farm located in northwest Oklahoma, the desired active power $P_d$ experiences large step changes between $0.3\operatorname{pu}$ and $1\operatorname{pu}$, and the desired reactive power $Q_d$ is such that the desired power factor $\text{PF}_d$ is fixed at $0.995$. As will be explained below, these values force the turbine to operate in both the MPT and PR modes, along with switching between them, under a realistic wind profile. Figure~\ref{fig:sim11} shows the simulation results for this scenario in both the per-unit and physical unit systems, wherever applicable. Subplot~1 shows the wind profile $V_w$. Subplot~2 shows the value of $C_p$, while subplot~3 shows the desired and actual active powers $P_d$ and $P$. Note that, for the first 1200 seconds during which $P_d$ is unachievable at $1\operatorname{pu}$, the turbine operates in the MPT mode and maximizes $P$, as indicated by the value of $C_p$ approaching its maximum of $0.48$ after a short transient (the turbine is initially at rest). At time 1200s when $P_d$ drops sharply from $1\operatorname{pu}$ to an achievable value of $0.3\operatorname{pu}$, the turbine quickly reduces the value of $C_p$, accurately regulates $P$ around $P_d$, and effectively rejects the ``disturbance'' $V_w$, thereby smoothly switches from the MPT mode to the PR mode. At time 2400s when $P_d$ goes from $0.3\operatorname{pu}$ back to $1\operatorname{pu}$, the MPT mode resumes. Because $V_w$ is strong enough at that time, $P$ approaches $P_d$. Moreover, the moment $P$ exceeds $P_d$ (which is equal to $P_{\text{rated}}$), the blade pitch angle $\beta$ increases as shown in subplot~8 in order to clip the power and protect the turbine. At time 3275s when $V_w$ becomes weaker, $\beta$ returns to $\beta_{\min}=0\operatorname{deg}$, thereby allowing the value of $C_p$ to return to its maximum of $0.48$ and $P$ to be maximized. Subplot~4 shows the desired and actual power factors $\text{PF}_d$ and $\text{PF}$, while subplots~5 and~6 show the desired and actual rotor angular velocities $\omega_{rd}$ and $\omega_r$ in normal and zoomed-in views. As can be seen from these subplots, throughout the simulation, both PF and $\omega_r$ are maintained near $\text{PF}_d$ and $\omega_{rd}$, respectively, affected only slightly by the random wind fluctuations. Moreover, the small S-shape variations in $\omega_{rd}$ in subplot~6 resemble those in Figure~\ref{fig:wrd}. Finally, subplots~7 and~8 show the control variables, i.e., the rotor voltages $v_{dr}$ and $v_{qr}$ as well as the blade pitch angle $\beta$.

The above simulation results suggest that the proposed controller not only is capable of operating effectively in both the MPT and PR modes, it is also capable of switching smoothly between them---all while not knowing the $C_p$-surface, air density, friction coefficient, and wind speed.

\section{Concluding Remarks}\label{sec:concl}

In this paper, we have designed a controller for a variable-speed wind turbine with a DFIG. The controller, consisting of four subcontrollers, has been developed based on a fifth-order, electromechanically-coupled, nonlinear model of the wind turbine by integrating several linear and nonlinear control strategies and exploiting time-scale separation in the dynamics. We have shown that the controller is able to make the wind turbine operate in both the MPT and PR modes and switch smoothly between them, in addition to maintaining a desired power factor. Furthermore, the controller does not require knowledge of the $C_p$-surface, air density, friction coefficient, and wind speed. Simulation has been carried out using a realistic wind profile, and the results demonstrate the capability and effectiveness of the controller. Future work includes designing a comprehensive wind farm controller for a multitude of turbines, which builds upon results documented in this paper.

\appendix[Proof of Theorem~\ref{thm:closedsys}]

First, we show (i). Setting $\dot{\omega}_r$ and $\dot{\hat{g}}$ in \eqref{eqn:closedsys1} and \eqref{eqn:closedsys2} to zero yields $g(\omega_r, \beta, V_w)=\max\{\hat{g}+\alpha\ln\frac{\omega_r}{\omega_{rd}}, 0\}$ and $\hat{g}=g(\omega_r, \beta, V_w)$. When $\hat{g}+\alpha\ln\frac{\omega_r}{\omega_{rd}}\geq0$, we have $\omega_r=\omega_{rd}$ and $\hat{g}=g(\omega_{rd}, \beta, V_w)$. Thus, $(\omega_{rd}, g(\omega_{rd}, \beta, V_w))$ is an equilibrium point, which is in $D$, since $0<\omega_{rd}\leq\omega_r^{(1)}$. On the other hand, when $\hat{g}+\alpha\ln\frac{\omega_r}{\omega_{rd}}<0$, we have $\omega_r\in\Omega$ and $\hat{g}=0$, where $\Omega=\{\omega\in(0,\infty):g(\omega, \beta, V_w)=0\}$ and $\omega_r^{(1)}=\min\Omega$. Since $\hat{g}+\alpha\ln\frac{\omega_r}{\omega_{rd}}<0$ and $\hat{g}=0$, we have $\omega_r<\omega_{rd}$. Since $\omega_r\in\Omega$, $\omega_r^{(1)}=\min\Omega$, and $\omega_{rd}\leq\omega_r^{(1)}$, we have $\omega_r\geq\omega_{rd}$. Hence, there is a contradiction, implying that when $\hat{g}+\alpha\ln\frac{\omega_r}{\omega_{rd}}<0$, there is no equilibrium point in $D$. This proves (i).

Next, we show (ii). To do so, it is useful to think of $D$ as a vertical strip in the two-dimensional state space $(\omega_r, \hat{g})$. Notice that on the right boundary of the strip where $\omega_r=\omega_r^{(1)}$, because of \eqref{eqn:closedsys1} and because $g(\omega_r^{(1)}, \beta, V_w)=0$ and $\max\{\hat{g}+\alpha\ln\frac{\omega_r^{(1)}}{\omega_{rd}},0\}\geq0$, we have $\dot{\omega}_r\leq0$. Thus, the state $(\omega_r, \hat{g})$ cannot escape $D$ through the right boundary. Next, note that for each fixed $\hat{g}\in\mathbb{R}$, there exists $\omega_r^{\star}>0$ such that for all $\omega_r\in(0, \omega_r^{\star})$, $\hat{g}+\alpha\ln\frac{\omega_r}{\omega_{rd}}<0$. This, along with \eqref{eqn:closedsys1} and the fact that $g(\omega, \beta, V_w)>0$ for all $\omega\in(0, \omega_r^{(1)})$, implies that near the left boundary of the strip where $\omega_r$ is arbitrarily small but positive, we have $\dot{\omega}_r>0$. Hence, the state $(\omega_r, \hat{g})$ cannot escape $D$ through the left boundary. This proves (ii).

Finally, we show (iii). Consider a Lyapunov function candidate $V:D\rightarrow\mathbb{R}$, defined as $V(\omega_r, \hat{g})=\alpha c(\omega_r\ln\frac{\omega_r}{\omega_{rd}}-\omega_r+\omega_{rd})+\frac{1}{2}(g(\omega_r, \beta, V_w)-\hat{g})^2$, where $c>0$ is to be determined. Note that $V$ is continuously differentiable over $D$. Moreover, $V$ is positive definite over $D$ with respect to the equilibrium point $(\omega_{rd}, g(\omega_{rd}, \beta, V_w))$, since $V(\omega_{rd}, g(\omega_{rd}, \beta, V_w))=0$ and $V(\omega_r, \hat{g})>0$ for all $(\omega_r, \hat{g})\neq(\omega_{rd}, g(\omega_{rd}, \beta, V_w))$ due to the property $\omega_r\ln\frac{\omega_r}{\omega_{rd}}-\omega_r+\omega_{rd}>0$ for all $\omega_r\neq\omega_{rd}$. Furthermore, $V$ is unbounded toward the top, bottom, and left boundary of the vertical strip $D$, but not so toward the right boundary of $D$. This is because for each fixed $\omega_r\in(0, \omega_r^{(1)}]$, $\lim_{|\hat{g}|\rightarrow\infty}V(\omega_r, \hat{g})=\infty$, and for each fixed $\hat{g}\in\mathbb{R}$, $\lim_{\omega_r\rightarrow0}V(\omega_r, \hat{g})=\infty$ and $V(\omega_r^{(1)}, \hat{g})<\infty$. Note that although $V$ is not unbounded toward the right boundary of $D$, the state $(\omega_r, \hat{g})$ cannot cross this boundary due to (ii).

Differentiating $V$ and using \eqref{eqn:closedsys1} and \eqref{eqn:closedsys2}, we get
\begin{align*}
J\dot{V}\!\!=\!\!\!\begin{bmatrix}
\alpha c\ln\frac{\omega_r}{\omega_{rd}}+(g-\hat{g})\frac{\partial g}{\partial \omega_r} \\ -(g-\hat{g})
\end{bmatrix}^T\!\!
\begin{bmatrix}
g-\max\{\hat{g}+\alpha\ln\frac{\omega_r}{\omega_{rd}},0\} \\ h(g-\hat{g})
\end{bmatrix}\!,
\end{align*}
where, for convenience, the function arguments are omitted. Note that because of (ii) and the above properties of $V$, to show (iii), it suffices to show that $\dot{V}$ is negative definite over $D$ with respect to the equilibrium point $(\omega_{rd}, g(\omega_{rd}, \beta, V_w))$. To this end, let $D$ be partitioned into two disjoint sets $D_1=\{(\omega_r, \hat{g})\in D: \hat{g}+\alpha\ln\frac{\omega_r}{\omega_{rd}}\geq0\}$ and $D_2=\{(\omega_r, \hat{g})\in D: \hat{g}+\alpha\ln\frac{\omega_r}{\omega_{rd}}<0\}$. Note that the equilibrium point $(\omega_{rd}, g(\omega_{rd}, \beta, V_w))$ is in $D_1$.

Suppose $(\omega_r, \hat{g})\in D_1$. Then, $\dot{V}$ takes a quadratic form:
\begin{align*}
J\dot{V}=
&=\!-\!\begin{bmatrix}
\ln\frac{\omega_r}{\omega_{rd}} \\ g-\hat{g}
\end{bmatrix}^T\!\!
\begin{bmatrix}
\alpha^2c & \frac{\alpha}{2}(\frac{\partial g}{\partial \omega_r}-c) \\
\frac{\alpha}{2}(\frac{\partial g}{\partial \omega_r}-c) & h-\frac{\partial g}{\partial \omega_r}
\end{bmatrix}\!
\begin{bmatrix}
\ln\frac{\omega_r}{\omega_{rd}} \\ g-\hat{g}
\end{bmatrix}.
\end{align*}
Note that if $(\omega_r, \hat{g})=(\omega_{rd}, g(\omega_{rd}, \beta, V_w))$, $\dot{V}=0$. Also, the leading principal minors of the above symmetric matrix are $\alpha^2c$ and $\alpha^2c(h-\frac{1}{4c}(\frac{\partial g}{\partial\omega_r}+c)^2)$. Thus, if $h$ and $c$ satisfy
\begin{align} 
h-\frac{1}{4c}(\frac{\partial}{\partial\omega}g(\omega,\beta,V_w)+c)^2>0, \quad\forall\omega\in(0, \infty),\label{eqn:case1}
\end{align}
then this symmetric matrix is positive definite, so that $\dot{V}<0$ for any $(\omega_r, \hat{g})\neq(\omega_{rd}, g(\omega_{rd}, \beta, V_w))$. Therefore, if $h$ and $c$ satisfy \eqref{eqn:case1}, $\dot{V}$ is negative definite over $D_1$ with respect to $(\omega_{rd}, g(\omega_{rd}, \beta, V_w))$.

Next, suppose $(\omega_r, \hat{g})\in D_2$. Then, $\dot{V}$ is bounded from above by a quadratic form:
\begin{align*}
J\dot{V}&=-h\hat{g}^2+(2h-\frac{\partial g}{\partial \omega_r})g\hat{g}+(\frac{\partial g}{\partial \omega_r}-h)g^2+\alpha cg\ln\frac{\omega_r}{\omega_{rd}}\displaybreak[0]\\
&\leq-h\hat{g}^2+(2h-\frac{\partial g}{\partial \omega_r})g\hat{g}+(\frac{\partial g}{\partial \omega_r}-h)g^2-cg\hat{g}\displaybreak[0]\\
&=-\begin{bmatrix}
\hat{g} \\ g
\end{bmatrix}^T
\begin{bmatrix}
h & \frac{1}{2}(\frac{\partial g}{\partial\omega_r}+c-2h) \\
\frac{1}{2}(\frac{\partial g}{\partial\omega_r}+c-2h) & h-\frac{\partial g}{\partial\omega_r}
\end{bmatrix}
\begin{bmatrix}
\hat{g} \\ g
\end{bmatrix}.
\end{align*}
Note that the leading principal minors of the above symmetric matrix are $h$ and $c(h-\frac{1}{4c}(\frac{\partial g}{\partial\omega_r}+c)^2)$. Thus, if $h$ and $c$ satisfy \eqref{eqn:case1}, then this symmetric matrix is positive definite. Since $(\omega_r, \hat{g})\in D_2$ and $\omega_{rd}\leq\omega_r^{(1)}$, if $\hat{g}=0$, then $g>0$. Thus, $\hat{g}$ and $g$ cannot be zero simultaneously. Hence, $\dot{V}<0$. Therefore, if $h$ and $c$ satisfy \eqref{eqn:case1}, $\dot{V}$ is negative over $D_2$.

As it follows from the above, if $h$ and $c$ satisfy \eqref{eqn:case1}, $\dot{V}$ is negative definite over $D$ with respect to the equilibrium point $(\omega_{rd}, g(\omega_{rd}, \beta, V_w))$, so that (iii) holds.

It remains to show that if $h$ satisfies \eqref{eq:condh}, then there exists $c>0$ such that \eqref{eqn:case1} holds. Suppose $h$ satisfies \eqref{eq:condh}. Let $F(\underline{\gamma}, \overline{\gamma})=\overline{\gamma}$ if $\overline{\gamma}\geq-\frac{1}{3}\underline{\gamma}$ and $F(\underline{\gamma}, \overline{\gamma})=-\frac{(\overline{\gamma}-\underline{\gamma})^2}{8(\underline{\gamma}+\overline{\gamma})}$ otherwise. Then, $h>F(\underline{\gamma}, \overline{\gamma})$. Let $f(x,\underline{\gamma},\overline{\gamma})=\frac{1}{4x}\max\{(\underline{\gamma}+x)^2, (\overline{\gamma}+x)^2\}$, where $x>0$. Then, it can be shown that $F(\underline{\gamma}, \overline{\gamma})=\min_{x>0}f(x,\underline{\gamma},\overline{\gamma})$ by considering the following three cases separately: $\overline{\gamma}\geq-\underline{\gamma}$, $-\underline{\gamma}>\overline{\gamma}\geq-\frac{1}{3}\underline{\gamma}$, and $-\frac{1}{3}\underline{\gamma}>\overline{\gamma}$. Because $h>F(\underline{\gamma}, \overline{\gamma})$, there exists $c>0$, given by $c=\arg\min_{x>0}f(x,\underline{\gamma},\overline{\gamma})$, such that $h>f(c,\underline{\gamma},\overline{\gamma})$. Because $\underline{\gamma}\leq\frac{\partial}{\partial\omega}g(\omega, \beta, V_w)\leq\overline{\gamma}$ for all $\omega\in(0,\infty)$ and by definition of $f(x,\underline{\gamma},\overline{\gamma})$, we have $\frac{1}{4c}(\frac{\partial}{\partial\omega}g(\omega, \beta, V_w)+c)^2\leq f(c,\underline{\gamma},\overline{\gamma})$ for all $\omega\in(0, \infty)$. Since $h>f(c,\underline{\gamma},\overline{\gamma})$, \eqref{eqn:case1} holds, as desired.

\bibliographystyle{IEEEtran}
\bibliography{paper_jrnl}

\end{document}